\documentclass[12pt,english,a4paper]{article}
\usepackage[T1]{fontenc}
\usepackage{amsmath}
\usepackage{graphicx}
\usepackage[symbol]{footmisc}
\usepackage{float}
\usepackage{amsmath,amsfonts}
\usepackage{mathrsfs}
\usepackage{amsmath}
\usepackage{authblk}

\allowdisplaybreaks

\def\R{\mathbb R}

\def\Liminf{\mathop{\underline{\lim}}\limits}

\def\Pb{\mathbf{P}}

\let\bar\overline
\let\tilde\widetilde
\def\limnto{\mathrel{\mathop{\longrightarrow\kern 0pt}\limits_{n\to\infty}}}

\let\bar\overline
\def\1{\mbox{1\hspace{-.25em}I}}

\def\Pb{\mathbf{P}}

\def\Ex{\mathbf{E}}

\def\BB{\mathbb{B}}

\def\KK{\mathbb{K}}

\def\RR{{\cal R}}

\def\II{\mathbb{I}}
\def\UU{\mathbb{U}}

\def\TT{\mathbb{T}}
 
\def\1{\mbox{1\hspace{-.25em}I}}

\let\tilde\widetilde

\newtheorem{theorem}{Theorem}
\newtheorem{lemma}{Lemma}

\newtheorem{definition}{Definition}

\usepackage{babel}

\begin{document}
\title{Poisson Source Localization on the Plane. Change-Point Case.}
\date{}

\author[1]{C. Farinetto}
\author[2]{Yu. A.  Kutoyants}
\author[3]{A. Top}
\affil[1,2,3]{\small  Le Mans University, Le Mans, France}
\affil[2]{Tomsk State  University, Tomsk, Russia,}
\affil[2]{National Research university, ''MPEI'',  Moscow, Russia}

\maketitle

\begin{abstract}
We present a detection problem where several spatially distributed sensors
observe Poisson signals emitted from a single  source of
unknown position. The
measurements at each sensor are modeled by independent inhomogeneous Poisson
processes. A method based on Bayesian change-point estimation is proposed to
identify the location of the source's coordinates. The asymptotic behavior of the
Bayesian estimator is studied. In particular the consistency and the asymptotic efficiency of the estimator are analyzed.
The limit distribution and the convergence of the moments are also described.
The similar statistical model could be used in GPS localization problems.
\end{abstract}

\bigskip{} \textbf{Key words}: Inhomogeneous Poisson process, change-point problem,
Bayesian estimator, likelihood ratio process,  source localization, sensors, GPS
localization.
 \bigskip{}
\date{}
\selectlanguage{english}

\maketitle

\bigskip{}

\section{Introduction}

In this work we study the properties of Bayesian estimators for the
localization of a  source emitting  Poisson  signals that propagate over an
area monitored by a set of sensors.  This mathematical model could be used for
the description of a radioactive emission, an explosion, a seismic activity or
the detection of weak optical signals. Sensors are electronic devices that can
measure changes in the environment around them, for instance there are light
sensors, proximity sensors, pressure sensors, heat sensors, radiation sensors
etc. The model under study could describe such data if the sequence of
observed random events is of Poisson nature.  Data obtained from a single
sensor is often not fully reliable and incomplete due to single device's
technical limitations. Using data from several sensors has advantages over
data collected from a single sensor. If several identical sensors are
employed, the observation process can be improved by combining individual
information to generate a more complete picture of the environment monitoring.
We refer the interested reader to Magee and Aggarwal \cite{MAGEE} or Chao
\cite{CHAO} for the advantages of using multiple sensors.  It has been shown
that the probability of measurement error decreases with the size of the
sensor network. However it is worth mentioning the complexity of the
monitoring system will increase with the number of sensors.  Source tracking
and localization is a problem of considerable importance that has attracted
the scientific interest.  Many examples of applications for such problems can
be found in environmental monitoring, industrial sensing, infrastructure
security, military tracking and diverse areas of security and defense, see for
instance Zhao \cite{ZG04} and Chong \cite{CK03}. 
 The present work focuses on
the detection of Poisson sources.  Due to the recent events security issues
have become more and more concerning and the problem of detecting radioactive
sources, more specifically the detection of illicit radioactive substances,
stored or in transit, has received great deal of attention by the engineering
community.

The detection of hidden nuclear material by means of sensors is an active area
of research as part of defensive strategies. One can consult the work of
Baidoo \cite{BAIDOO}, Liu \cite{LIU} and Rao \cite{RAO} for details and
references on this topic.  Nuclear radiations are a probabilistic physical
process consisting of discrete emissions of particles that can be recorded by
radiation sensors. Those emissions have been mathematically modeled with help
of Poisson point processes which provide natural models describing their
properties, see for instance Evans \cite{EVANS} or Knoll \cite{KNOLL}.  Apart
from radiation measurements, typical examples on the use of Poisson point
processes include modeling streams of photo-electrons produced by light on
photosensitive surfaces \cite{MANDEL}, laser radar detection and ranging of
objects \cite{KAR}, earthquake aftershocks \cite{OGATA}, electrical response
of nerves to stimulus \cite{SNYDER} and others, for application to tracking
and sensing we refer to the book of Streit \cite{STREIT}.  Special cases of
the source localization problem have been studied in the past, for instance
Howse \cite{HOW} described least squares estimation algorithms to estimate the
location of a possibly moving source by a fixed number of sensors. For
multiple sources Maximum Likelihood Estimation (MLE) was considered by
Morelande \cite{MORELANDE}.  An iterative procedure for calculating MLEs of a
single nuclear source from radiation measurements as well as corresponding
Cramer-Rao bounds for localization accuracy was given by Baidoo \cite{BAIDOO}.
Concerning Bayesian statistics Liu \cite{LIU} presented a technique to locate
a source according to Bayesian update methods. The results of Pahlajani
\cite{PAHLAJANI} are also noteworthy: their paper studies the presence of a
source using Likelihood ratio calculation and a Neyman-Pearson test. In what
follows we suppose that there is a single source generating a signal. Our goal
is to describe the asymptotic behavior of the Bayesian estimator (BE) of its
coordinates through the method developed by Ibragimov and Khasminskii
\cite{IH81} for the study of such estimators. We show that the rate of
convergence of the estimator is $n$ and that the limit distribution is not
Gaussian. A lower bound on the mean-square risk is proposed and the BE is
proved to be asymptotically efficient.

Note that the same mathematical model can be used in the problem of
GPS-localization on the plane \cite{XL13}. Indeed, in this case the signals
are emitted by $k$ fixed emitters and an object receiving these signals has to
define its own position. Here once more we have $k$ signals with unknown
moments of arriving and using the estimators of these moments the object can
construct the estimator of its position.

\section{ Statement of the  problem}

We are interested in locating the source of an event with the help of several
spatially distributed independent sensors monitoring an area over a fixed
time-interval.  For example, if we have a radioactive source, then each sensor
records ambient measurements, for instance radiations due to natural isotopes
in the environment. When the event occurs, then the sensors record the sum of
ambient measurements and the measurements related to the event. The two
signals are independently, and we consider that each sensor records a single
inhomogeneous Poisson process whose intensity is the sum of the intensities
due to both ambient and background event measurements.

Popular network topologies for source localization problems that were
considered in other studies are grids of sensors \cite{LIU} and triangular
arrays \cite{CHIN}.  In order to identify the source location we use a configuration of
sensors forming a triangle.


  In our case, we have  sequences of
measurements from three sensors and collected within the same time window. The
measurements from each sensor are sent to a central processing unit (fusion
center) that combines the data and estimates the coordinates of the source.


\bigskip

\begin{figure}[ht]
\hspace{3cm}\includegraphics[width=9cm]   {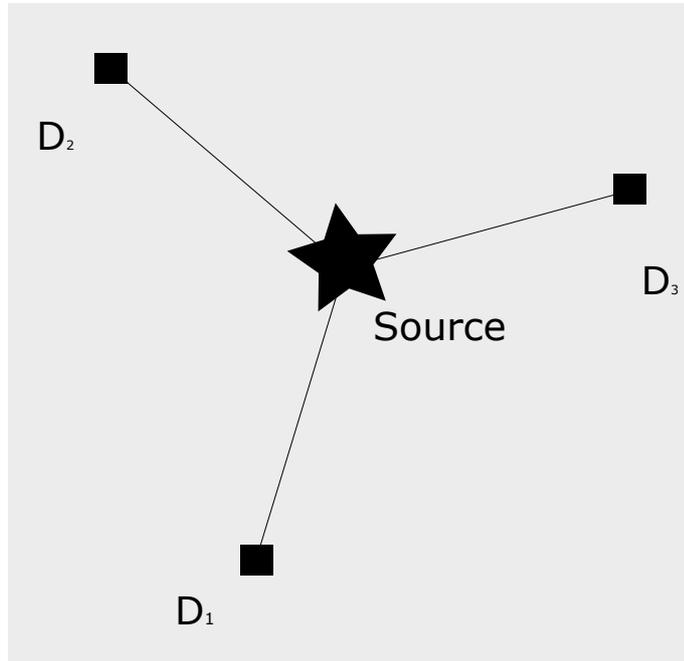}
\caption{Model of observations}
\end{figure}

The source is located at an unknown position $D_0$ with coordinates
$\vartheta_0= (x_0,y_0) $ inside a convex set $\Theta \subset {\cal R}^2$.
Three sensors are placed in the field at known positions at points
$D_1,D_2,D_3$ with the coordinates $\vartheta_j= (x_j,y_j)$, $j=1, 2, 3$.
Each sensor records on the time interval $[0,T]$ a signal modeled by a Poisson
point process $X_{j}\!=\!\bigl\{ X_{j}(t),\,0\leq t \leq T \bigr\}$, $j= 1, 2,
3$ of intensity function $\lambda_j \left(\vartheta_0, t\right),0\leq t\leq
T$.  These intensity functions are supposed to be of the form
\begin{align*}
\lambda_j \left(\vartheta_0, t\right)=\lambda \left(t-\tau _j\right)+\lambda
_0,\qquad 0\leq t\leq T. 
\end{align*}
Here $\lambda _0>0$ is a known intensity of the background noise, $\lambda
\left(t\right)$ is the known intensity function of the signal and $\tau
_j=\tau _j\left(\vartheta _0\right)$ is the arrival time of the signal to the
$j$-th sensor (delay). This delay is calculated following the usual rule
\begin{align}
\label{tau}
  \tau_j\left(\vartheta _0\right) = \frac{||\vartheta_j - \vartheta_0||}{\nu},
\end{align}
where $||\cdot||$ is the Euclidean norm and $\nu$ is the known rate of
propagation of the signal in the monitored area.  We suppose that $\lambda
\left(t\right)=0$ for $t<0$. At time $t=0$ the emission of signals begins and
$\tau _j$ is the arrival time of the signal to the $j$-th sensor. We are
concerned by estimating the position $\vartheta _0$ of the Poisson
source. We are interested in the models of observations which allow the
estimation with small errors such that $\Ex_{\vartheta_0 }\left(\bar \vartheta
-\vartheta _0\right)^2=o\left(1\right) $. As usual such situations are
considered in an asymptotic framework. The small errors can be obtained if the
intensity of the signal takes large values or a periodical Poisson process
could describe the data. Another possibility is to have many sensors.  We take
the model with large intensity functions $\lambda_j \left(\vartheta_0,
t\right)=\lambda_{j,n} \left(\vartheta_0, t\right) $ which can be written as
follows
\begin{align}
\label{int}
\lambda_{j,n} \left(\vartheta_0, t\right)=n\lambda \left(t-\tau _j\right)+n\lambda
_0,\qquad 0\leq t\leq T. 
\end{align}
Here $n$ is a ``large parameter'' and we study estimators as $n\rightarrow
\infty $. For example, such a model could be obtained in the case of three clusters, where
each cluster includes $n$ detectors.

The likelihood ratio function $L\left(\vartheta ,X^n\right) $  is
\begin{align}
\label{lr}
\ln L\left(\vartheta ,X^n\right)&=\sum_{j=1}^{3}\int_{\tau _j}^{T}\ln \left(1+\frac{\lambda
  \left(t-\tau _j\right)}{\lambda _0}\right) {\rm
  d}X_j\left(t\right) -n\sum_{j=1}^{3}\int_{\tau _j}^{T}\lambda
  \left(t-\tau _j\right){\rm d}t. 
\end{align}
Here $\tau _j=\tau _j\left(\vartheta \right)$ and $X^n=\left(X_j\left(t\right),0\leq
t\leq T, j=1,2,3\right)$ are counting processes from three detectors. Based on this
likelihood ratio formula we define the  maximum likelihood estimator (MLE)
$\hat\vartheta _n$   and Bayesian estimator (BE)  
$\tilde\vartheta _n$ by 
\begin{align}
\label{mle}
L\left(\hat\vartheta _n ,X^n\right)=\sup_{\vartheta \in \Theta
}L\left(\vartheta ,X^n\right),
\end{align}
and
\begin{align}
\label{be}
\tilde\vartheta _n=\frac{\int_{\Theta }^{}\vartheta p\left(\vartheta
  \right)L\left(\vartheta ,X^n\right){\rm d}\vartheta }{\int_{\Theta }^{}
  p\left(\vartheta \right)L\left(\vartheta ,X^n\right){\rm d}\vartheta },
\end{align}
respectively.

Here $p\left(\vartheta \right),\vartheta \in \Theta $ is the prior density. As
the limit properties of the BE do not depend on the prior density, we could
consider a non-informative prior such as the uniform density. For any other
positive continuous function $p\left(\cdot \right)$ the limit properties will
remain the same.

Recall that in the case of a discontinuous intensity function $\lambda
\left(\cdot \right)$ the definition of the MLE has to be modified since
\begin{align*}
\ln L\left(\vartheta ,X^n\right)&=\sum_{j=1}^{3}\sum_{i=1}^{N_j}\ln
\left(1+\frac{\lambda \left(t_{i,j}-\tau _j\left(\vartheta
  \right)\right)}{\lambda _0}\right)\\
&\qquad \qquad  -n\sum_{j=1}^{3}\int_{\tau
  _j\left(\vartheta \right)}^{T}\lambda \left(t-\tau _j\left(\vartheta
\right)\right){\rm d}t.
\end{align*}
Here $t_{i,j},i=1,\ldots,N_j$ are the registration times of the events in
the $j$-th sensor and $N_j$ is the total number of events in this
sensor. Of course, if $N_j=0$, then we set
\begin{align*}
\sum_{i=1}^{N_j}\ln \left(1+\frac{\lambda
  \left(t_{i,j}-\tau _j\left(\vartheta \right)\right)}{\lambda _0}\right)=0.
\end{align*}
We write formally 
\begin{align*}
\max \left(L(\hat\vartheta _n- ,X^n),L(\hat\vartheta _n+
,X^n) \right)=\sup_{\vartheta \in \Theta }L(\vartheta ,X^n) 
\end{align*}
which we understand as follows. The function 
$$
M\left(\tau _1\left(\vartheta \right),\tau _2\left(\vartheta
\right),\tau _3\left(\vartheta \right),X^n \right)=L\left(\vartheta
,X^n\right)
$$ 
has jumps at points   $\tau _j\left(\vartheta \right)=t_{i,j} $ and its
supremum is in one of the jump points. It can be written 
\begin{align*}
\sup_{\vartheta \in\Theta } L\left(\vartheta
,X^n\right)=\max M\left(\tau _1(\hat\vartheta_n )\pm,\tau _2(\hat\vartheta_n 
)\pm,\tau _3(\hat\vartheta_n )\pm,X^n  \right).
\end{align*}
 Here $M\left(\tau _1(\vartheta )\pm,\tau _2(\vartheta
)\pm,\tau _2(\vartheta )\pm,X^n  \right)$ are left and
right limits  of the function $M\left(\tau _1(\vartheta ),\tau _2(\vartheta
),\tau _2(\vartheta ),X^n  \right)$ at the points $\tau _j(\vartheta ) $.

There are several different types of problems associated with the
identification of the location depending on the regularity of the function
$\lambda \left(t \right)$. In particular, the rate of convergence of the mean
square error of the estimators $\bar\vartheta _n$ is
\begin{align*}
\Ex_{\vartheta_0 }\left(\bar
\vartheta_n -\vartheta _0\right)^2=\frac{C}{n^\gamma }\left(1+o\left(1\right)\right),
\end{align*}
where the parameter $\gamma >0$ depends on the regularity of the function
$\lambda \left(\cdot \right)$.

Let us present three of them. All the cases are illustrated using the following model
\begin{align}
\label{if}
\lambda\left(\vartheta ,t\right)=2\left|\frac{t-\tau _j\left(\vartheta \right)
}{\delta }\right|^\kappa \1_{\left\{0\leq t-\tau _j\left(\vartheta
  \right)<\delta \right\}} +2\1_{\left\{ t-\tau _j\left(\vartheta \right)\geq
  \delta \right\}} +1.
\end{align}
Statistical problems related to different types of regularity could be
obtained according to the values of parameter $\kappa$.

\bigskip

{\bf \textbullet Smooth case.} Suppose that the function $\lambda \left( \cdot
\right) $ 
in \eqref{int} is  continuously differentiable, then the problem of parameter
estimation is regular.

The MLE $\hat\vartheta _n$ and BE $\tilde\vartheta _n$ (under regularity
conditions) are consistent, asymptotically normal
\begin{align*}
\sqrt{n}\left(\hat\vartheta _n-\vartheta _0 \right)\Longrightarrow {\cal
  N}\left(0,\II\left(\vartheta _0\right)^{-1}\right),\quad
\sqrt{n}\left(\tilde\vartheta _n-\vartheta _0 \right)\Longrightarrow {\cal 
  N}\left(0,\II\left(\vartheta _0\right)^{-1}\right),
\end{align*}
the moments converge and   both estimators are asymptotically
efficient. Here $\II\left(\vartheta _0\right) $ represents the Fisher information
matrix. For the mean square error the following relation holds true:
\begin{align*} 
\Ex_{\vartheta_0 }\left\|\hat
\vartheta_n -\vartheta _0\right\|^2=\frac{C}{n }\left(1+o\left(1\right)\right),
\end{align*}
i.e., $\gamma =1$.  

This case corresponds to the intensity function \eqref{if} with
$\kappa>\frac{1}{2}$. An example of such an intensity function is given in
Fig. 2.


\begin{figure}[h]
\hspace{3cm}\includegraphics[width=9cm]   {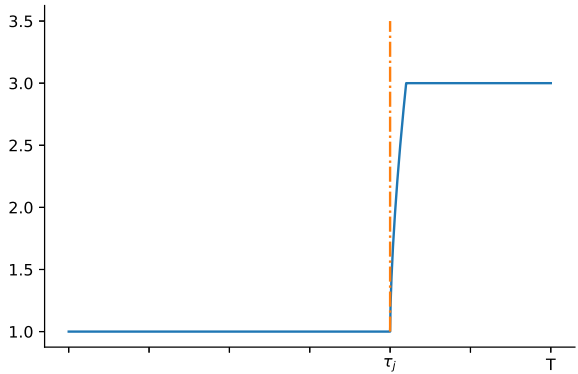}
\caption{Intensity (6) with $\kappa = \frac{5}{8}, \delta =0,1$.}
\end{figure}

It is worth mentioning that the derivative of this function is a discontinuous
function; however it is continuous in $L_2\left(0,T\right)$ and the MLE has
all the aforementioned properties.

We describe these properties of estimators in the problem
of Poisson source localization   in the forthcoming work \cite{CK18}.

\bigskip

{\bf \textbullet Change-point  case.} Suppose that the intensity function in \eqref{int} has the
  following form
\begin{align*}
\lambda_{j,n} \left(\vartheta ,t\right)=n\lambda_1\left(t-\tau_j
\right)\1_{\left\{t\geq \tau _j\right\}} +n\lambda  _0,\qquad 0\leq t\leq T. 
\end{align*}
Here $\lambda _1\left(t\right)>0, t\geq 0$ and $\lambda  _0>0$
are known.

This type of statistical problems corresponds to the intensity function
\eqref{if} with $\kappa =0$ and $\delta =0$ (see Fig. 3).

{\bf Here is Fig.3}

\begin{figure}[h]
\hspace{3cm}\includegraphics[width=9cm]   {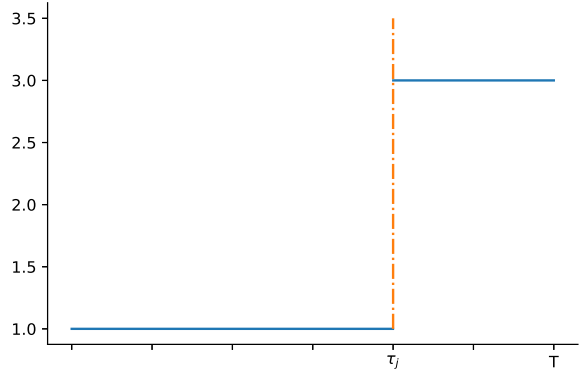}
\caption{Intensity (6) with $ \kappa =0, \delta =0$.}
\end{figure}

 In this situation the intensities of the observed Poisson
 processes have positive jumps equal to $n\lambda _1\left(0\right)$ at the points
 $t=\tau _j=\tau _j\left(\vartheta_0 \right)$. This is a non regular parameter
 estimation problem, where the MLE and BE  have the
 normalization $n$ and different limit distributions
\begin{align*}
{n}\left(\hat\vartheta _n-\vartheta _0 \right)\Longrightarrow \hat \zeta ,\qquad
{n}\left(\tilde\vartheta _n-\vartheta _0 \right)\Longrightarrow \tilde\zeta .
\end{align*}
The moments of these estimators converge, but only the BE is asymptotically
efficient. The random vectors $\hat\zeta $ and $\tilde\zeta $ are exponential
functionals of some Poisson processes.  The mean square error decreases as
follows
\begin{align*}
\Ex_{\vartheta_0 } \left\|\tilde
\vartheta_n -\vartheta _0\right\| ^2=\frac{C}{n^2 }\left(1+o\left(1\right)\right),
\end{align*}
i.e., $\gamma =2$. Similar results in the case of an  one-dimensional parameter
$\vartheta $ could be found in \cite{Kut98}.

Here we focus on the study of the BE for this model of observations.

\bigskip

{\bf \textbullet Cusp-type  case.} This case is in some sense intermediate between the
  smooth and change-point cases. Suppose that the intensity function has the
  following form
\begin{align*}
\lambda_{j,n} \left(\vartheta,t\right)=n\,\lambda _1\left(t-\tau
_j\right)\left|\frac{t-\tau_j 
}{\delta } 
\right|^\kappa\1_{\left\{0\leq t-\tau_j\leq  \delta \right\}} 
+n\lambda _1\left(t- \tau_j\right)\1_{\left\{t- \tau_j> \delta \right\}}+n\lambda _0
 .
\end{align*}
The parameter $\kappa \in (0,\frac{1}{2})$, the parameter $\delta >0$ takes
small values and the function $\lambda _1\left(t \right)>0$.

An example of such a function is given in Fig.4.
 
{\bf Here is Fig.4}

\begin{figure}[ht]
\hspace{3cm}\includegraphics[width=9cm]   {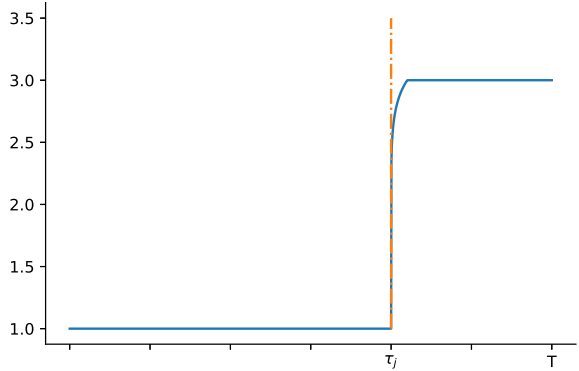}
\caption{Intensity (6) with $\kappa = 0,1, \delta =0,1.$}
\end{figure}
\bigskip 

In the statistical literature change-point  models are well studied, but in
some real cases  the intensity function could have
jumps since due to the physical laws the electrical current can not have
jumps and the cusp-type model fits much better to the real data with strongly
increased intensities. The intensity of the signal increases from zero to
$\lambda \left(\tau +\delta \right)$ in the small interval $\left[\tau ,\tau
  +\delta \right]$.   Note that for these values of $\kappa $ the Fisher
information does not exist which leads to a singular estimation problem. The MLE and
BE for this model of observations are consistent, have different limit
distributions
\begin{align*}
{n}^{\frac{1}{2\kappa +1}}\left(\hat\vartheta _n-\vartheta _0
\right)\Longrightarrow \hat \xi  ,\qquad 
{n}^{\frac{1}{2\kappa +1}}\left(\tilde\vartheta _n-\vartheta _0
\right)\Longrightarrow \tilde\zeta ,
\end{align*}
  the moments  converge and only the BE is asymptotically efficient. The
  random vectors  $\hat \xi  $  and $\tilde\xi  $ are  
exponential functionals of the fractional Brownian motions.  
\begin{align*}
\Ex_{\vartheta_0 } \left\|\hat
\vartheta_n -\vartheta _0 \right\|^2=\frac{C}{n^{\frac{2}{2\kappa +1}}
}\left(1+o\left(1\right)\right), 
\end{align*}
i.e., $\gamma =\frac{2}{2\kappa +1} $  and   $1<\gamma <2$. These cases 
will be
studied in the forthcoming  work \cite{CDK18}. For the one-dimensional parameter case see
\cite{D03}. 

The properties of the MLE and BE of the one-dimensional parameter in such
three types of regularity problems for the signals observed in the white
Gaussian noise are discussed in  \cite{DKKN18}.

 \section{Main results}
There are three sensors with coordinates $\vartheta_j=\left(x_j,y_j\right) ,
j=1,2,3$ which measure the particles emitted by some source at the point
$\vartheta _0=\left(x_0,y_0\right) $. 
The observations are modeled by three independent  inhomogeneous
Poisson processes 
$X^n=\left(X_j\left(t\right),0\leq t\leq T, j=1,2,3\right)$ with respective intensity
functions
\begin{align*}
\lambda _{j,n}\left(\vartheta_0 ,t\right)=n\lambda\left(t-\tau
_j\left(\vartheta \right)\right)\1_{\left\{t\geq \tau 
  _j\left(\vartheta \right)\right\}}+n\lambda _0 ,\qquad 0\leq t\leq T
\end{align*}
where $\lambda \left(t\right)>0$ and $\lambda _0>0$. The arrival times
of the signals in the $j$-th sensor according to \eqref{tau} are $\tau
_j=\tau _j\left(\vartheta _0\right)$ and the position of
the source $\vartheta_0 =\left(x_0,y_0\right) \in\Theta \subset {\RR}^2$ will be estimated.  We
suppose that $\Theta =\left(\alpha _1,\alpha _2\right)\times \left(\beta
_1,\beta _2\right)$ with finite $\alpha _i,\beta _i$. The set $\Theta $ is
bounded, open and convex.  Of course, we suppose that for all $\vartheta \in
\Theta $ the corresponding $\tau \left(\vartheta \right)\in \left(0,T\right)$.

Note that if the model of observations with the constant intensities
of the signal and noise is considered, i.e.,
\begin{align*}
\lambda _{j,n}\left(\vartheta_0 ,t\right)=n\lambda_1\1_{\left\{t\geq \tau 
  _j\left(\vartheta \right)\right\}}+n\lambda _0 ,\qquad 0\leq t\leq T,
\end{align*}
where $\lambda _1=\lambda \left(0\right)>0$,  then the asymptotic properties
of the estimators  will be the same.

We study the asymptotic ($n\rightarrow \infty$) behavior of the Bayesian
estimator of the unknown parameter $\vartheta_0= (x_0,y_0) $.  It is worth noticing that in such non
regular estimation problems the asymptotic 
results could be applied even for moderate values of $n$ since we have faster
convergence of estimators  (rate $n$ and not $\sqrt{n}$ as in the  regular case). 

Let us introduce the quantities
\begin{align*}
\underline{\tau }=\min_{j=1,2,3}\inf_{\vartheta \in\Theta }\tau
_j\left(\vartheta \right),\qquad  \bar \tau =\max_{j=1,2,3}\sup_{\vartheta \in\Theta }\tau
_j\left(\vartheta \right),\qquad \TT=\left[0,T-\bar\tau \right].
\end{align*}
At this point we have to suppose some conditions providing the identifiabilily of the position of
the source.

Conditions ${\cal I}$:
\begin{description}
\item[${\cal I}_1$.] {\it The location of the source is different from the sensor location.
   Consequently we suppose that 
  there exists a small constant $\varepsilon >0$ such that for every possible
  position of the source $\vartheta_0\in\Theta $ and} $j=1, 2,3$
 $$ 
\rho_j= || \vartheta_j -   \vartheta_0 || \geq \varepsilon.
 $$

\item[${\cal I}_2$.] {\it The function $\lambda \left(s\right), s\in \TT$ has
  two continuous derivatives}

\item[${\cal I}_3$.] {\it The sensors are not aligned, therefore}
\begin{equation*}
\begin{vmatrix}
x_1 & x_2 &  x_3 \\
y_1 & y_2 &  y_3 \\
1 & 1 & 1
\end{vmatrix} \neq  0.
\end{equation*}
\end{description}
By condition ${\cal I}_1$ the case $\tau _j=0$ is excluded. 
Due to  condition ${\cal I}_1$  we restrain the parameter space to
$$\Theta= \Big [ (\alpha_1,\alpha_2)\times(\beta_1,\beta_2) \Big ] \backslash
\Bigg [ \bigcup_{j=1}^{3} B(\vartheta_j, \varepsilon ) \Bigg ], $$ where
$B(\vartheta_j,\varepsilon) = \{ z \in \mathbb{R}^2: ||\vartheta_j -z|| \leq
\varepsilon \}$. If the position of the source coincides with the position of
one of the sensors, then for this sensor $\tau _j=0$ and the properties of the
estimators will be different. For example, the limit likelihood ratio
$Z\left(u\right)$ can be defined for the positive values of one component of
$u$ only. This situation corresponds to the case, where the true value of the
unknown parameter is on the border of a parametric set (see, e.g. \cite{Kut98},
where such situation was described).  Remark, that if the condition ${\cal
  I}_2 $ is not fulfilled and the sensors are in the same line, then the
consistent estimation of the position $\vartheta _0$ is not feasible. Of course
such conclusion depends on the set $\Theta $ too. Suppose that the detectors
are on a line on the seashore and the source can be only be located on one side, then
two detectors are sufficient for the consistent estimation of the position of
the Poisson (radioactive) source.

The likelihood  $L\left( \vartheta, X^n\right) $ according to \eqref{lr}  is
given by (see for example \cite{Kut98}). 
\begin{align*}
&\ln L\left( \vartheta, X^{(n)}\right)=\sum _{j=1}^{3} \int_{0}^{T}\ln
\frac{\lambda_{j,n }\left(\vartheta, t\right)}{n\lambda _0}{\rm d}X_j(t) -\sum
_{j=1}^{3} \int_{0}^{T} \left(\lambda_{j,n} \left(\vartheta, t\right)-n\lambda
_0\right){\rm d}t \\
 &\qquad \qquad = \sum _{j=1}^{3} \int_{\tau _j}^{T}\ln \bigl ( 1 +
\frac{\lambda\left(t-\tau _j\right)}{\lambda _0} \bigr ){\rm
  d}X_j(t) -n\sum _{j=1}^{3} \int_{\tau _j}^{T} \lambda\left(t-\tau _j\right){\rm d}t .
\end{align*}
Recall that here $\tau _j=\tau _j\left(\vartheta  \right)$.

If the intensity function of the signal is constant $\lambda
\left(t\right)\equiv \lambda _1>0$, then the likelihood ratio is simplified
\begin{align*}
\ln L\left( \vartheta, X^{(n)}\right)=\ln\left(1+\frac{\lambda_1}{\lambda
  _0}\right)\sum _{j=1}^{3} \left[X_j\left(T\right)-X_j\left(\tau
  _j\right) \right]-n\lambda _1 \sum
_{j=1}^{3}\left[T-\tau _j\right].
\end{align*}

The Bayesian estimator
$\tilde\vartheta_n=(\tilde x_n, \tilde y_n) $ of the parameter $\vartheta_0=(
x_0,y_0) $ with respect to the quadratic loss function is defined by a
conditional expectation which can be written as follows
\begin{align*}
\tilde\vartheta_n= \mathbf{E}\left(\vartheta /
X^{(n)}\right)=\int_{  \Theta } \vartheta p(\vartheta)L\left( \vartheta,
X^{(n)}\right)\;{\rm d}\vartheta \Biggl ( \int_ {\Theta} p(\vartheta)L\left(
\vartheta, X^{(n)}\right)\;{\rm d}\vartheta\Biggr)^{-1}.
\end{align*}

Even if the vector $\vartheta $ is not random with a given prior  density we can
 use this formula to calculate $\tilde \vartheta _n$ which
is no more a conditional expectation, but just some way to construct the
estimator. In this case it can be called {\it generalized Bayesian
  estimator} \cite{IH81}. Therefore we can take  any positive
continuous function $p\left(\vartheta  \right),\vartheta \in\Theta $. 
For example, as the set $\Theta $ is bounded, we can put $p\left(\vartheta
\right)=1$. 

Note that if the intensity of the signal is constant $\lambda
\left(t\right)\equiv \lambda _1$, then the estimator can be calculated as follows
\begin{align*}
\tilde\vartheta _n= \frac{\int_{\Theta }^{}\vartheta
  \prod_{j=1}^{3}\bigl(1+\frac{\lambda _1}{\lambda _0}\bigr)^{-X_j\left(\tau
    _j\left(\vartheta \right)\right)} e^{n\lambda _1\sum_{j=1}^{3}\tau _j\left(\vartheta \right)}{\rm
    d}\vartheta}{\int_{\Theta }^{} \prod_{j=1}^{3}\bigl(1+\frac{\lambda
    _1}{\lambda _0}\bigr)^{-X_j\left(\tau _j\left(\vartheta \right)\right)} e^{n\lambda
    _1\sum_{j=1}^{3}\tau _j\left(\vartheta \right)}{\rm d}\vartheta },
\end{align*}
where $\tau _j\left(\vartheta \right)=\nu ^{-1}\left\|\vartheta_j-\vartheta
\right\|$. 

In order to describe the properties of the Bayesian estimator, we need some
additional notations. First let us introduce the unit vectors $m_j$, for $j=1,\cdots,3$
\begin{equation*}
m_j= \left( \frac{x_j-x_0}{\rho_j},\frac{y_j-y_0}{\rho_j}\right),\qquad
\rho _j=\left\|\vartheta_j -\vartheta _0\right\|,\quad \left\|m_j \right\|=1
\end{equation*}
and the sets
\begin{align*}
\BB_j=\left\{u:\quad \langle  m_j,u \rangle \geq 0
\right\},\quad \BB_j^c=\left\{u:\quad \langle  m_j,u \rangle <
0  \right\} .
\end{align*}
Here $\langle m_j,u \rangle$ denotes the Euclidean scalar product of the
vectors $m_j$ and $u=\left(u_1,u_2\right) $.
The limit likelihood ratio $Z\left(u\right),u\in \R^2$ we denote as follows
\begin{align*}
\ln Z\left(u\right)&=\ell \sum_{j=1}^{3}\left[ \Pi
_{j,+}\left(u\right)\1_{\left\{u\in\BB_j\right\}} -\Pi
_{j,-}\left(u\right)\1_{\left\{u\in\BB_j^c\right\}}   \right]\\
&\qquad \qquad 
-\frac{\lambda _1 }{\nu }\langle  m_1+m_2+m_3,u \rangle ,
\end{align*}
where $\ell=\ln \left(1+\frac{\lambda _1}{\lambda _0}\right)$, $\Pi
_{j,+}\left(u\right),u\in \BB_j $ and $\Pi _{j,-}\left(u\right),u\in \BB_j^c$
are independent Poisson random fields such that
\begin{align*}
\Ex_{\vartheta _0}\Pi _{j,+}\left(u\right)=\frac{\langle  m_j,u \rangle}{\nu },\qquad
\Ex_{\vartheta _0}\Pi _{j,+}\left(u\right)=-\frac{\langle  m_j,u \rangle }{\nu }.
\end{align*}

Second, we define the random vector $\tilde \zeta=(\tilde\zeta_1,\tilde\zeta_2) $ with the components
\begin{equation*}
 \tilde \zeta_1 =  \int_{{\cal R}^2}u_1 Z(u_1,u_2) \;{\rm d} u_1 {\rm d} u_2
\Biggl (\int \int_{{\cal R}^2} Z(u_1,u_2) \;{\rm d} u_1  {\rm d} u_2\Biggr)^{-1}
\end{equation*}
and
\begin{equation*}
 \tilde \zeta_2 =  \int_{{\cal R}^2}u_2 Z(u_1,u_2) \; {\rm d} u_1 {\rm d} u_2
\Biggl (\int \int_{{\cal R}^2} Z(u_1,u_2) \;{\rm d} u_1 {\rm d} u_2\Biggr)^{-1}.
\end{equation*}

The main results of this work are 
the following two theorems.   We first introduce the lower
bound on the risk of all estimators.
\begin{theorem}
\label{T1}
Let the conditions ${\cal I}$ be fulfilled. Then for all $\vartheta_0 \in
\Theta$ and a quadratic loss function,
\begin{eqnarray}
\label{lb}
\lim_{\delta \rightarrow 0}\Liminf_{n\rightarrow \infty }
\inf_{\bar{\vartheta}_n}\sup_{||\vartheta-\vartheta_0||<\delta}
n^{2}\mathbf{E}_{\vartheta}  \left\|\bar{\vartheta}_n-\vartheta\right\|^2 \geq
\mathbf{E} \|\tilde\zeta\|^2 .  
\end{eqnarray}
\end{theorem}
Here  the $inf$ is taken over all
possible estimators $\bar \vartheta_n$ of the
parameter $\vartheta$.  The inequality (\ref{lb}) allows us to give the
following definition of  efficient estimator.
\begin{definition}
\textit{Let the conditions} ${\cal I}$ \textit{be satisfied. The
  estimator $\vartheta^{*}_n$ is asymptotically efficient, if for all
  $\vartheta_0 \in \Theta$ we have}
\begin{align}
\lim_{\delta \rightarrow 0}\underset{n \rightarrow
  +\infty}{\lim}\sup_{||\vartheta -\vartheta_0||<\delta}
n^{2}\mathbf{E}_{\vartheta} \left\|\vartheta^{*}_n-\vartheta\right\|^2 =\mathbf{E}
 \|\tilde\zeta \|^2.
 \label{D1}
\end{align}
\end{definition}

The second theorem describes the asymptotic behavior of the estimator
$\tilde\vartheta_n=(\tilde x_n, \tilde y_n)$.
 \begin{theorem}
\label{T2}
Let the conditions ${\cal I}$ be fulfilled. Then the Bayesian estimator
$\tilde{\vartheta}_n$ is uniformly on compacts $\KK\subset\Theta $ consistent:
for any $\gamma  >0$
\begin{align*}
\sup_{\vartheta _0\in \KK} \Pb_{\vartheta
  _0}\left(\|\tilde{\vartheta}_n-\vartheta_0 \|>\gamma
\right)\longrightarrow 0,  
\end{align*}
we have  convergence in distribution
\begin{equation*} 
n\left(\tilde{\vartheta}_n-\vartheta_0\right)\quad \Longrightarrow \quad  
 \tilde\zeta,
\end{equation*}
and convergence of  moments: for any $p>0$
\begin{align*}
\lim_{n\rightarrow \infty }n^p \Ex_{\vartheta
  _0}\|\tilde{\vartheta}_n-\vartheta_0\|^p =\Ex_{\vartheta
  _0}\|\;\tilde\zeta \;\|^p ,
\end{align*}
and $\tilde{\vartheta}_n$ is  asymptotically efficient.
\end{theorem}
The proofs of these theorems are given in the next section. 
They are  based on the general results of Ibragimov and
Khasminskii  \cite{IH81} for the problem of parameter estimation in the case
of i.i.d. observations with a discontinuous density function and  the
application of their results to the study of  Bayesian estimators
for inhomogeneous Poisson processes see (\cite{Kut98}, Chapter 5). 

Let us remind the main steps of these proofs. Introduce the normalized
likelihood ratio random field
\begin{align*}
Z_n\left(u\right)=\frac{L\left(\vartheta _0+\frac{u}{n},
  X^n\right)}{L\left(\vartheta _0, X^n\right)},\qquad u\in \UU_n,
\end{align*}
where 
\begin{align*}
\UU_n=\left\{u:\;\vartheta _0+\frac{u}{n}\in\Theta \right\}.
\end{align*}
Moreover we extend the set $\UU_n$ to cover the balls around the sensors
\begin{align*}
\UU_n=\Bigl(n\left(\alpha _1-x_0\right),n\left(\alpha _2-x_0\right)\Bigr)\times
\Bigl( n\left(\beta _1-y_0\right),n\left(\beta_2-y_0\right) \Bigr)\nearrow \RR^2,
\end{align*}
as $n\rightarrow \infty $. i.e., we extended the process $Z_n\left(u\right)$
on the values $u$ belonging 
to the balls $\vartheta _0+\frac{u}{n} \in B\left(\vartheta_j,\varepsilon
\right)$. This requires certain modifications of the general method developed
in \cite{IH81},  which can
be done without difficulties. What is important is to respect the condition
$\vartheta _0\not\in B\left(\vartheta_j,\varepsilon
\right)$. 

Suppose that we have already proved the convergence of finite dimensional
distributions $Z_n\left(\cdot \right)\Longrightarrow Z\left(\cdot
\right)$. Below we change the variables $\vartheta =\vartheta
_0+\frac{u}{n}$. We have
\begin{align*}
\tilde\vartheta _n&={\int_{\Theta }^{}\vartheta \frac{L\left(\vartheta
    ,X^n\right)}{L\left(\vartheta_0
    ,X^n\right) }\;{\rm d}\vartheta }\left({\int_{\Theta }^{} \frac{L\left(\vartheta
    ,X^n\right)}{L\left(\vartheta_0
    ,X^n\right) }\;{\rm d}\vartheta}\right)^{-1}\\
&=\vartheta _0+\frac{1}{n}\int_{\UU_n}^{}uZ_n\left(u\right)\,{\rm
  d}u\left(\int_{\UU_n}^{}Z_n\left(u\right)\,{\rm d}u \right)^{-1} 
\end{align*}
and
\begin{align*}
n\left(\tilde\vartheta _n-\vartheta
_0\right)=\int_{\UU_n}^{}uZ_n\left(u\right)\,{\rm
  d}u\left(\int_{\UU_n}^{}Z_n\left(u\right)\,{\rm d}u \right)^{-1}. 
\end{align*}
If we prove the convergence
\begin{align*}
&\left(\int_{\UU_n}^{}u_1Z_n\left(u\right){\rm
  d}u,\int_{\UU_n}^{}u_2Z_n\left(u\right){\rm
  d}u, \int_{\UU_n}^{} Z_n\left(u\right)\,{\rm 
  d}u\right)\\
&\qquad \qquad \qquad \qquad \Longrightarrow \left(\int_{\RR^2}^{}u_1Z\left(u\right){\rm
    d}u,\int_{\RR^2}^{}u_2Z\left(u\right){\rm   d}u,\int_{\RR^2}^{} Z\left(u\right)\,{\rm   d}u \right),
\end{align*}
then we obtain the limit 
\begin{align*}
n\left(\tilde\vartheta _n-\vartheta _0\right)\Longrightarrow \tilde\zeta .
\end{align*}
To obtain the convergence of moments we have to check the uniform
integrability of the random variables $\left\|n\left(\tilde\vartheta
_n-\vartheta _0\right)\right\|^p $ for any $p>0$. 

This work was realized in \cite{IH81}  in a  sufficiently general framework
(see Theorem 1.10.2 there). In the next section we verify the conditions
of this theorem. 

Suppose that we already proved  Theorem \ref{T2}, then the proof of 
Theorem \ref{T1} could be done as follows. Let us fix some small $\delta >0$, then
\begin{align*}
\sup_{||\vartheta-\vartheta_0||<\delta}
n^{2}\mathbf{E}_{\vartheta}\|\bar{\vartheta}_n-\vartheta\|^2&\geq
n^{2}\int_{B\left(\vartheta _0,\delta
  \right)}^{}\mathbf{E}_{\vartheta}\|\bar{\vartheta}_n-\vartheta\|^2
q\left(\vartheta  
\right){\rm d}\vartheta\\ 
&\geq n^{2}\int_{B\left(\vartheta _0,\delta
  \right)}^{}\mathbf{E}_{\vartheta}
\|\tilde{\vartheta}_{q,n}-\vartheta\|^2q\left(\vartheta
\right){\rm d}\vartheta,
\end{align*}
where $q\left(\vartheta \right),\vartheta \in B\left(\vartheta _0,\delta
\right)$ is some positive continuous density on $B\left(\vartheta _0,\delta
\right) $  and $\tilde{\vartheta}_{q,n} $ is a BE, which
corresponds to this prior density. From the convergence of second  moments  we
have
\begin{align*} 
n^{2}\int_{B\left(\vartheta _0,\delta
  \right)}^{}\mathbf{E}_{\vartheta}
\|\tilde{\vartheta}_{q,n}-\vartheta \|^2q\left(\vartheta
\right){\rm d}\vartheta\longrightarrow \int_{B\left(\vartheta _0,\delta
  \right)}^{}\mathbf{E}_{\vartheta}
\|\tilde\zeta \|^2\;q\left(\vartheta
\right)\,{\rm d}\vartheta.
\end{align*}
The continuity of $\mathbf{E}_{\vartheta} \|\tilde\zeta \|^2 $
w.r.t. $\vartheta $ allows us to write the last limit
\begin{align*}
\int_{B\left(\vartheta _0,\delta
  \right)}^{}\mathbf{E}_{\vartheta}
\|\tilde\zeta  \|^2\;q\left(\vartheta
\right)\,{\rm d}\vartheta\longrightarrow \mathbf{E}_{\vartheta_0}
\|\tilde\zeta  \|^2
\end{align*}
as $\delta \rightarrow 0$. Note that the lower bound \eqref{lb} is a
particular case of more general result in \cite{IH81}.

\section{Proofs}

Introduce the normalized likelihood random field
\begin{eqnarray*}
Z_{n}(u) &=&
\exp \Biggl\{ \sum _{j=1}^{3} \int_{0}^{T}\ln \frac{ \lambda_{j,n}
  \left(\vartheta_0+\frac{u}{n}, t\right) }{\lambda_{j,n} \left(\vartheta_0,
  t\right) }{\rm d} X_j(t)\\ &&- \sum _{j=1}^{3} \int_{0}^{T} \left(
\lambda_{j,n}(\vartheta_0+\frac{u}{n}, t)-\lambda_{j,n} \left(\vartheta_0,
t\right)\right) {\rm d}t \Biggr\},
\end{eqnarray*}
where  $u=(u_1,u_2) \in  \UU_n$.
\begin{lemma}
\label{L1}
Let the conditions  ${\cal  I}_1,{\cal  I}_2$ be satisfied, then the finite dimensional
distributions of the process $Z_{ n }(u),u\in\UU_n$ converge to the finite
dimensional distributions of the process $Z(u),u\in\RR^2$
 and this convergence is uniform with respect to $\vartheta_0 \in {\KK}$.
\end{lemma}
{\bf Proof.}
The characteristic function of $\ln Z_{n}(u)$ is calculated as follows (see
\cite{Kut98}) 
\begin{align*}
&\Phi_n(\mu;u)=\mathbf{E}_{\vartheta_0} \exp \left [i \mu \ln Z_n(u) \right
  ]\\ &\qquad =\exp \biggl\{ \sum _{j=1}^{3}\int_{0}^{T}\biggl [ \exp \biggl(i
    \mu \ln \frac{\lambda_{j,n} \left(\vartheta_0+\frac{u}{n},
      t\right)}{\lambda_ {j,n} \left(\vartheta_0, t\right)} \biggr)-1\biggr
  ]\lambda_{j,n} \left(\vartheta_0, t\right) {\rm d}t\\ &\qquad \quad -i \mu
  \sum _{j=1}^{3}\int_{0}^{T}\left( \lambda_ {j,n}(\vartheta_0+\frac{u}{n},
  t)-\lambda_ {j,n} \left(\vartheta_0, t\right)\right) {\rm d}t\biggr\}.
\end{align*}
Introduce  the  sets $A_k^n$ for $k=1,\cdots,8$, and $u=(u_1,u_2) \in \UU_n $
\begin{align*}
 A^n_1&= \{ u  \in \UU_n , \quad\langle u, m_1\rangle \geq 0 ,
\langle u, m_2\rangle \leq 0, \quad \langle u, m_3\rangle \leq 0\},\\
 A^n_2&= \{ u  \in \UU_n , \quad \langle u, m_1\rangle \geq 0
,\langle u, m_2\rangle \geq 0, \quad \langle u, m_3\rangle \leq 0\},\\
 A^n_3&= \{ u  \in \UU_n , \quad \langle u, m_1\rangle \geq 0 ,
\langle u, m_2\rangle \geq 0, \quad\langle u, m_3\rangle \geq 0\},\\
 A^n_4&= \{ u  \in \UU_n , \quad \langle u, m_1\rangle \leq 0
,\langle u, m_2\rangle \geq 0, \quad \langle u, m_3\rangle \geq 0\},\\
 A^n_5&= \{ u  \in \UU_n , \quad \langle u, m_1\rangle \leq0 ,
\langle u, m_2\rangle \leq 0, \quad \langle u, m_3\rangle \geq 0\},\\
 A^n_6&= \{ u \in \UU_n , \quad \langle u, m_1\rangle \leq 0 ,
\langle u, m_2\rangle \leq 0, \quad \langle u, m_3\rangle \leq 0\}.\\
 A^n_7&= \{ u  \in \UU_n, \quad \langle u, m_1\rangle \geq 0 ,
\langle u, m_2\rangle < 0, \quad \langle u, m_3\rangle \geq 0 \}.\\
 A^n_8&= \{ u  \in \UU_n, \quad \langle u, m_1\rangle < 0 ,\langle u,
m_2\rangle \geq 0, \quad \langle u, m_3\rangle < 0\}.
\end{align*} 
Define $\vartheta _u=\vartheta _0+\frac{u}{n}$, $\tau _j=\tau
_j\left(\vartheta _0\right)$, $ \rho _j=\nu \tau _j$ and 
\begin{eqnarray*}
\tau _j(\vartheta _u)=\frac{1}{\nu }\sqrt{\left (x_j-x_0 - \frac{u_1}{n}
  \right)^2 + \left(y_j-y_0 - 
  \frac{u_2}{n} \right )^2 }.
\end{eqnarray*}

It follows from condition ${\cal I}_1$ that $\tau _j(\vartheta _u)$ is
differentiable w.r.t. $u$ on
$\UU_n$. Using the Taylor expansion we obtain
\begin{eqnarray*}
\tau _j(\vartheta _u)&=&\tau _j-\frac{ u_1(x_j-x_0) + u_2(y_j-y_0) }{\nu  n \rho_j}
+\varepsilon_n(u)\\
 &=&\tau _j-\frac{1}{\nu  n}\langle u,m_j \rangle
+\varepsilon_n(u),
\end{eqnarray*}
where  $n\varepsilon_n(u) \rightarrow 0 $
uniformly on compacts $u$ as $n\rightarrow \infty$.
Thus
\begin{equation*}
  \tau _j(\vartheta _u)-\tau _j=-\frac{1}{\nu  n}\langle u,m_j \rangle +\varepsilon_n(u).
\end{equation*}
Therefore for all $j=1, 2, 3$, bounded sets of $u$ and $n$ sufficiently large
we have
\begin{equation*}
\left \{ \begin{array}{ll} \tau _j \geq \tau _j(\vartheta _u)  , \quad  if \quad
  \langle u,m_j \rangle \geq 0, \\ 
 \tau _j \leq \tau _j(\vartheta _u)  , \quad if \quad \langle u,m_j \rangle \leq 0.
\end{array} \right. \end{equation*}
We will use this fact to calculate the characteristic function $\Phi_n(\mu;u)$
for each set $A^n_k$ , $k=1,\cdots,8$ and obtain its limit.

 If $ u \in A^n_1$, then $\tau _1\geq \tau _1\left(\vartheta _u\right)$, $\tau
 _2\leq \tau _2\left(\vartheta _u\right) $ and  $\tau
 _3\leq \tau _3\left(\vartheta _u\right) $. Therefore we can write 
\begin{align*}
&\int_{0}^{T}\biggl [ \exp \biggl(i
    \mu \ln \frac{\lambda_{1,n} \left(\vartheta_0+\frac{u}{n},
      t\right)}{\lambda_ {1,n} \left(\vartheta_0, t\right)} \biggr)-1\biggr
  ]\lambda_{1,n} \left(\vartheta_0, t\right) {\rm d}t\\
&\qquad =n\lambda _0\int_{\tau _1\left(\vartheta _u\right)}^{\tau _1}\biggl [ \exp \biggl(i
    \mu \ln \frac{\lambda \left(t-\tau _1\left(\vartheta _u\right)\right)+\lambda _0}{\lambda _0} \biggr)-1\biggr
  ] {\rm d}t\\
&\qquad \quad +n\int_{\tau _1}^{T}\biggl [ \exp \biggl(i
    \mu \ln \frac{\lambda \left(t-\tau _1\left(\vartheta _u\right)\right)+\lambda _0}{\lambda \left(t-\tau _1\right)+\lambda _0} \biggr)-1\biggr
  ]\left[\lambda \left(t-\tau _1\right)+\lambda _0\right] {\rm d}t.
\end{align*}
Using once again Taylor's expansions by the powers of $\frac{u}{n}$ we obtain
the representation
\begin{align*}
&\int_{0}^{T}\biggl [ \exp \biggl(i
    \mu \ln \frac{\lambda_{1,n} \left(\vartheta_0+\frac{u}{n},
      t\right)}{\lambda_ {1,n} \left(\vartheta_0, t\right)} \biggr)-1\biggr
  ]\lambda_{1,n} \left(\vartheta_0, t\right) {\rm d}t\\
&\qquad \qquad = \left[\exp\left\{i\mu \ln\frac{\lambda _1+\lambda _0}{\lambda
    _0}\right\}-1\right]\frac{\lambda _0}{\nu }\langle u, m_1 \rangle +o\left(1\right).
\end{align*}
The similar arguments give us the relations
\begin{align*}
&\int_{0}^{T}\biggl [ \exp \biggl(i
    \mu \ln \frac{\lambda_{2,n} \left(\vartheta_0+\frac{u}{n},
      t\right)}{\lambda_ {2,n} \left(\vartheta_0, t\right)} \biggr)-1\biggr
  ]\lambda_{2,n} \left(\vartheta_0, t\right) {\rm d}t\\
&\qquad \qquad = -\left[\exp\left\{-i\mu \ln\frac{\lambda _1+\lambda _0}{\lambda
    _0}\right\}-1\right]\frac{\lambda _1+\lambda _0}{\nu }\langle u, m_2 \rangle +o\left(1\right)
\end{align*}
and
\begin{align*}
&\int_{0}^{T}\biggl [ \exp \biggl(i
    \mu \ln \frac{\lambda_{3,n} \left(\vartheta_0+\frac{u}{n},
      t\right)}{\lambda_ {3,n} \left(\vartheta_0, t\right)} \biggr)-1\biggr
  ]\lambda_{3,n} \left(\vartheta_0, t\right) {\rm d}t\\
&\qquad \qquad = -\left[\exp\left\{-i\mu \ln\frac{\lambda _1+\lambda _0}{\lambda
    _0}\right\}-1\right]\frac{\lambda _1+\lambda _0}{\nu }\langle u, m_3
  \rangle +o\left(1\right). 
\end{align*}

Therefore for $u\in A_1^n$ we  obtain the limit
\begin{align*}
& \lim_{n \rightarrow \infty} \Phi_n(\mu;u) = \exp
\biggl\{\biggl[ \exp \big(i \mu \ell  \big )-1\biggr
]\frac{\lambda_0}{\nu}\langle u, m_1 \rangle\\
&\qquad \qquad  -\biggl[ \exp \big(-i\mu\ell  \big)-1\biggr
]\frac{\lambda_0+\lambda_1}{\nu } \langle u, m_2+m_3 \rangle -i \mu r(u)
\biggr\}. 
\end{align*}

If $ u \in A^n_2$, then  similar arguments allow us to verify that 
\begin{align*}
& \lim_{n \rightarrow \infty} \Phi_n(\mu;u) = \exp
\biggl\{\biggl[ \exp \big(i \mu \ell  \big ) -1\biggr
]\frac{\lambda_0}{\nu}\langle u, m_1 +m_2\rangle \\
&\qquad \qquad  -\biggl[ \exp \big(-i \mu \ell \big)-1\biggr
]\frac{\lambda_0+\lambda_1}{\nu} \langle u, m_3 \rangle -i \mu r(u)
\biggr\}.
\end{align*}
For  $ u \in A^n_3$ we have
 \begin{align*}
\lim_{n \rightarrow \infty} \Phi_n(\mu;u) = \exp
\biggl\{\biggl[ \exp \big(i \mu \ell  \big)-1\biggr
]\frac{\lambda_0}{\nu}\langle u, m_1 +m_2+m_3\rangle -i \mu r(u) \biggr\}.
\end{align*}
For  other sets $A_k^n$ we have the corresponding limits. For all sets
these limits provide  the convergence of characteristic functions
\begin{align*}
\Ex_{\vartheta _0}\exp \left[i \mu \ln Z_n\left(u\right)\right]\quad
\longrightarrow \quad 
\Ex_{\vartheta _0}\exp \left[i \mu \ln Z \left(u\right)\right].
\end{align*}
Therefore we have the convergence  of one-dimensional distributions.  

Using the same arguments it is possible to verify the convergence of the
finite-dimensional distributions too, i.e., for any $u_1,\ldots,u_L$ and 
reals $\mu _1,\ldots,\mu _L$ we have
\begin{align*}
\Ex_{\vartheta _0}\exp \left[i\sum_{l=1}^{L} \mu_l \ln Z_n\left(u_l\right)\right]\quad
\longrightarrow \quad 
\Ex_{\vartheta _0}\exp \left[i \sum_{l=1}^{L}\mu_l \ln Z \left(u_l\right)\right].
\end{align*}
Moreover from the presented proofs it follows that the convergence of
finite-dimensional distributions is uniform on the compacts $\KK\subset \Theta
$. In particular,
\begin{align*}
\lim_{n\rightarrow \infty }\sup_{\vartheta _0\in\KK}\left|\Ex_{\vartheta _0}\exp \left[i\sum_{l=1}^{L} \mu_l \ln Z_n\left(u_l\right)\right]-
\Ex_{\vartheta _0}\exp \left[i \sum_{l=1}^{L}\mu_l \ln Z
  \left(u_l\right)\right]  \right|=0. 
\end{align*}

Further we need  the following result.
\begin{lemma}
\label{L2}
Let the  condition ${\cal  I}_1,{\cal  I}_2$ be fulfilled, then for any $R>0$ and $\left\|u\right\|+\left\|
v\right\|\leq R,\;  u,v\in \UU_n$   we have 
\begin{align*}
\sup_{\vartheta _0\in\KK}{\bf E}_{\vartheta_0} \left | Z_{n}^{\frac{1}{2}}(u) -
Z_{n}^{\frac{1}{2}}(v) \, \right | ^{2} \leq
C\left(1+R\right)\left\|u-v\right\|,
\end{align*}
 where $C>0.$
\end{lemma}
{\bf Proof.}
According to the Lemma 1.5 in  \cite{Kut98}, we have
\begin{align*}
&{\bf E}_{\vartheta_0} \left | Z_{n}^{\frac{1}{2}}(u) - Z_{n}^{\frac{1}{2}}(v)
  \, \right | ^{2} \\ 
&\quad \leq \sum _{j=1}^{3}\int_{0}^{T}\biggl [
    \sqrt{\lambda_{j,n} \left(\vartheta_0+\frac{u}{n}, t\right)}
    -\sqrt{\lambda_{j,n} \left(\vartheta_0+\frac{v}{n}, t\right)} \biggr]^2
     {\rm d}t \\
 &\quad \leq n \sum _{j=1}^{3}\int_{0}^{T} \biggl [
       \sqrt{ \lambda\left(t-\tau _j\left(\vartheta _u\right)\right) \1_{
           \left\{ t >\tau _j(\vartheta _u))            \right\} }+\lambda_0  }\\
&\quad \qquad  \qquad \qquad \qquad \qquad  - \sqrt{\lambda\left(t-\tau _j\left(\vartheta _v\right)\right) \1_{ \left\{ t >
           \tau _j(\vartheta _v) \right\} }+\lambda_0  } \biggr]^2 {\rm d}t\\
 &\quad
     \leq Cn \sum _{j=1}^{3}\int_{0}^{T} {  \left [\lambda\left(t-\tau
         _j\left(\vartheta _u\right)\right)\1_{ 
           \left\{ t >\tau _j(\vartheta _u) \right\} } - \lambda\left(t-\tau
         _j\left(\vartheta _v\right)\right)\1_{ \left\{ t > 
           \tau _j(\vartheta _v) \right\} } \right ]^2 } {\rm d}t.
\end{align*}
Here we use the elementary relations
\begin{align*}
\left[\sqrt{a}-\sqrt{b}\right]^2=\frac{\left[a-b\right]^2}{\left[\sqrt{a}
+\sqrt{b}\right]^2}\leq C\,\left[a-b\right]^2,\qquad C=\frac{1}{4 M},
\end{align*}
where  $a>0,b>0$ and  $M\geq \left(a\vee b\right)$.

Consider the values $\left|u\right|+\left|v\right|\leq R$ with some
$R>0$. Then using once again Taylor's expansions we obtain 
\begin{align*}
\lambda\left(t-\tau_j\left(\vartheta
_u\right)\right)-\lambda\left(t-\tau_j\left(\vartheta
_v\right)\right)= \frac{1}{\nu n}  \lambda'\left(t-\tau_j\right) \langle u-v,m_j\rangle+\varepsilon _n\left(u,v\right)
\end{align*}
and for large $n$
\begin{align*}
\left|\tau_j\left(\vartheta _u\right)-\tau _j\left(\vartheta _v\right)\right|\leq \frac{2}{\nu
  n}\left|\langle u-v,m_j\rangle\right|\leq \frac{C}{n}\left\|u-v\right\|.
\end{align*}
These two estimates allow us to write 
\begin{align*}
&n\sum_{j=1}^{3}\int_{0}^{T} {  \left [\lambda\left(t-\tau
         _j\left(\vartheta _u\right)\right)\1_{ 
           \left\{ t >\tau _j(\vartheta _u) \right\} } - \lambda\left(t-\tau
         _j\left(\vartheta _v\right)\right)\1_{ \left\{ t > 
           \tau _j(\vartheta _v) \right\} } \right ]^2 } {\rm d}t\\
&\qquad \qquad \leq C\left\|u-v\right\|+\frac{C}{n}\left\|u-v\right\|^2\leq
  C\left(1+R\right)\left\|u-v\right\| .
\end{align*}

The last result is given in  the next lemma.
\begin{lemma}
\label{L3}
Let conditions $ {\cal I}$  be fulfilled, then for $u \in \UU_n$
\begin{align}
\label{10}
\sup_{\vartheta _0\in\KK}{\bf E}_{\vartheta_0}   Z_{n}^{\frac{1}{2}}(u) \leq
 e^{-\kappa  ||u||  },
\end{align}
where $\kappa >0$.
\end{lemma}
{\bf Proof.} According to Lemma 1.5 of  \cite{Kut98} we can write
\begin{equation*}
{\bf E}_{\vartheta_0} \Big [ Z_{n}^{\frac{1}{2}}(u) \Big ] = \exp \Biggl \{ -
\frac{1}{2} \sum _{j=1}^{3}\int_{0}^{T}\biggl [ \sqrt{\lambda_{j,n}
    \left(\vartheta_0+\frac{u}{n}, t\right)} -\sqrt{\lambda_{j,n}
    \left(\vartheta_0, t\right)} \biggr]^2 {\rm d}t \Biggr \}.
 \end{equation*}
Elementary  calculations leads to 
\begin{align*}
 &\biggl [ \sqrt{\lambda_{j,n }\left(\vartheta_0+\frac{u}{n}, t\right)}
    -\sqrt{\lambda_{j,n }\left(\vartheta_0, t\right)}\biggr]^2 \\
 &\qquad
  =\frac{ n \left [\lambda \left(t-\tau _j\left(\vartheta _u\right)\right)\1_{
        \left\{ t > \vartheta _j(\vartheta _u) 
        \right\} } -\lambda \left(t-\tau _j\right) \1_{ \left\{ t > \tau _j \right\} } \right
    ]^2 }{\bigg [ \sqrt{ \lambda \left(t-\tau _j\left(\vartheta
        _u\right)\right) \1_{ \left\{ t > \tau _j(\vartheta _u) \right\}
        }+\lambda_0  } + \sqrt{ \lambda \left(t-\tau _j\right) \1_{ \left\{ t 
          > \tau _j \right\} }+\lambda_0  } \bigg ]^2}\\
 &\qquad \geq n c \left [\lambda \left(t-\tau _j\left(\vartheta _u\right)\right)\1_{
        \left\{ t > \vartheta _j(\vartheta _u) 
        \right\} } -\lambda \left(t-\tau _j\right) \1_{ \left\{ t > \tau _j \right\} } \right
    ]^2
\end{align*}
where $c =  \frac{ 1}{4 \lambda _M}$  with the constant $\lambda _M\geq
\lambda \left(t\right)+\lambda _0$. 

Let us now consider $\vartheta $ such that $\left\|\vartheta -\vartheta
_0\right\|\leq \delta $ with small $\delta >0$ and such that $\tau
_j\left(\vartheta \right)>\tau _j$. Then  for sufficiently small $\delta $ we
can write
\begin{align*}
&\int_{0}^{T}\left [\lambda \left(t-\tau _j\left(\vartheta \right)\right)\1_{
        \left\{ t > \tau _j(\vartheta ) 
        \right\} } -\lambda \left(t-\tau _j\right) \1_{ \left\{ t > \tau _j \right\} } \right
    ]^2{\rm d}t \\
&\qquad \qquad =\int_{\tau _j}^{\tau _j\left(\vartheta \right) }\lambda
  \left(t-\tau _j\right) ^2{\rm d}t+ \int_{\tau _j\left(\vartheta \right)
  }^{T}\left[\lambda \left(t-\tau _j\left(\vartheta \right)\right) -\lambda
    \left(t-\tau _j\right)  \right]^2{\rm d}t \\
&\qquad \qquad \geq k \left[\tau _j\left(\vartheta \right)-\tau _j
    \right]-c_j\left\|\vartheta -\vartheta _0\right\| ^2\geq k_j\left[\tau _j\left(\vartheta \right)-\tau _j
    \right]
\end{align*}
with $k =\min_t\lambda \left(t\right)^2>0$ and  some positive constant  $ c_j,k_j$.

Using this last inequality we obtain
\begin{align}
\label{8}
&\sum _{j=1}^{3}\int_{0}^{T} \biggl[ \sqrt{\lambda_{j,n}
    \left(\vartheta_0+\frac{u}{n}, t\right)} -\sqrt{\lambda_{j,n}
    \left(\vartheta_0, t\right)} \biggr]^{2} {\rm d}t\geq n \gamma \sum
_{j=1}^{3} \left | \tau _j(\vartheta ) - \tau _j \left(\vartheta
_0\right)\right |\nonumber\\ 
&\qquad \geq \gamma \sum _{j=1}^{3} \left
| \langle m_j, u \rangle \right | + \varepsilon_n(\delta )\geq  \gamma_1 \sum _{j=1}^{3} \left |
\Big < m_j, \frac{u}{\left\|u\right\|} \Big > \right | \left\|u\right\|\nonumber\\ 
&\qquad  \geq\gamma_1 \inf_{\left\|e\right\|=1}\sum _{j=1}^{3} \left | 
\langle m_j, e \rangle\right | ||u||\geq
 \kappa _1 \left\|u\right\|  ,
\end{align}
where $\kappa _1  >0$.

Next we consider the case $\left\|\vartheta -\vartheta
_0\right\|=\left\|\frac{u}{n}\right\| > \delta$. Let us denote
\begin{align*}
g\left(\vartheta _0,\delta \right)=\inf_{\left\|\vartheta -\vartheta
  _0\right\|>\delta } \sum _{j=1}^{3}\int_{0}^{T}\left [\lambda \left(t-\tau
  _j\left(\vartheta \right)\right)\1_{ 
        \left\{ t > \tau _j(\vartheta ) 
        \right\} } -\lambda \left(t-\tau _j\right) \1_{ \left\{ t > \tau _j \right\} } \right
    ]^2{\rm d}t.
\end{align*}
Remark that for any compact $\KK\subset \Theta $
\begin{align*}
g_\KK\left(\delta \right)=\inf_{\vartheta _0\in\KK}g\left(\vartheta _0,\delta \right)>0.
\end{align*}
Indeed, if $g_\KK\left(\delta \right)=0 $, then there exists $\vartheta _1\not=
\vartheta _0$, such that 
\begin{align*}
\sum _{j=1}^{3}\int_{0}^{T}\left [\lambda \left(t-\tau _j\left(\vartheta_1
  \right)\right)\1_{ \left\{ t > \tau _j(\vartheta_1 ) \right\} } -\lambda
  \left(t-\tau _j\left(\vartheta _0\right)\right) \1_{ \left\{ t > \tau
    _j\left(\vartheta _0\right) \right\} } \right ]^2{\rm d}t=0 .
\end{align*}
Due to the indicator functions this equality is possible iff $\tau _j\left(\vartheta
_1\right)=\tau _j\left(\vartheta
_0\right),j=1,2,3 $ but from the geometrical consideration this  is
impossible. Therefore $g_\KK\left(\delta \right)>0 $ and for $\left\|\vartheta
-\vartheta _0\right\|\geq \delta $ we can write
\begin{align}
\label{9}
&n\sum _{j=1}^{3}\int_{0}^{T}\left [\lambda \left(t-\tau _j\left(\vartheta
  \right)\right)\1_{ \left\{ t > \tau _j(\vartheta ) \right\} } -\lambda
  \left(t-\tau _j\left(\vartheta _0\right)\right) \1_{ \left\{ t > \tau
    _j\left(\vartheta _0\right) \right\} } \right ]^2{\rm d}t\nonumber\\
&\qquad \qquad \geq ng_\KK\left(\delta \right)\geq ng_\KK\left(\delta
  \right)\frac{\left\|\vartheta -\vartheta _0\right\|}{D\left(\Theta \right)}
  \geq \kappa _2\left\|u\right\|. 
\end{align}
Here 
$$
D\left(\Theta \right)=\sup_{\vartheta ,\vartheta _0\in\Theta }\left\|\vartheta
-\vartheta _0\right\|,\qquad  \kappa _2=\frac{g_\KK\left(\delta \right)}{D\left(\Theta \right)}.
$$
From the estimates \eqref{8} and \eqref{9} it follows that there exists
$\kappa >0$ such that
\begin{align*}
\sum _{j=1}^{3}\int_{0}^{T}\biggl [ \sqrt{\lambda_{j,n}
    \left(\vartheta_0+\frac{u}{n}, t\right)} -\sqrt{\lambda_{j,n}
    \left(\vartheta_0, t\right)} \biggr]^2 {\rm d}t\geq 2\kappa \left\|u\right\|.
\end{align*}

This last estimate proves \eqref{10}. 

The properties of the normalized likelihood ratio
$Z_n\left(u\right),u\in\UU_n$ described in the Lemmas \ref{L1}-\ref{L3} allow
us to cite Theorem 1.10.2 in \cite{IH81} and according to this theorem the BE
$\tilde\vartheta _n$ has all the  properties mentioned in Theorem \ref{T1}.

\section{Simulations}

We illustrate the convergence of the estimators  by means of numerical simulations.
Consider the problem of localization of a Poisson source at the point
$\vartheta_0=(0,0)$. We have   three sensors $\vartheta_j$ $(j =1,2,3)$
respectively located at coordinates 
$\vartheta _1=(8.5,0) $, $\vartheta _2=(0,8.5) $ and $\vartheta _3=\big (8.5
\cos(\frac{5 \pi}{4}), 8.5 \sin \big (\frac{5 
  \pi}{4}) \big ) $. We choose the values $\lambda_0 =1$, $\lambda_1=2$ and
for convenience $\nu = 1$.  Each sensor located at position $\vartheta_j$
records in the fixed time interval $[0,10]$ measurements that are modeled by a Poisson
point processes of intensity function
$$
\lambda_j \left(\vartheta_0, t\right) = n + 2 n \1_{\left\{t \geq \tau_j
  \right\}}.
$$ 
The parameter space of the unknown coordinates of the source
$\vartheta_0$ was chosen as $\Theta =\left( -1,1\right)\times
\left(-1,1\right)$ and the prior density of $\vartheta_0$ is the uniform density in  the unit square, i.e. $p\left(\vartheta \right) = \frac{1}{4}
\1_{\left\{(x,y) \in [-1,1]^2 \right\}} $. The
BE $\tilde\vartheta _{n}$ was calculated using simulations for $n$ running in
the range $[1,100]$.  Fig. 5  displays the evolution of the Euclidean
distance between the BE $\tilde{ \vartheta}_{n} =(\tilde{x}_n,
\tilde{y}_{n}) $ and $\vartheta_0$ with $n$. 


\begin{figure}[ht]
\hspace{3cm}\includegraphics[width=10cm, height = 10cm]   {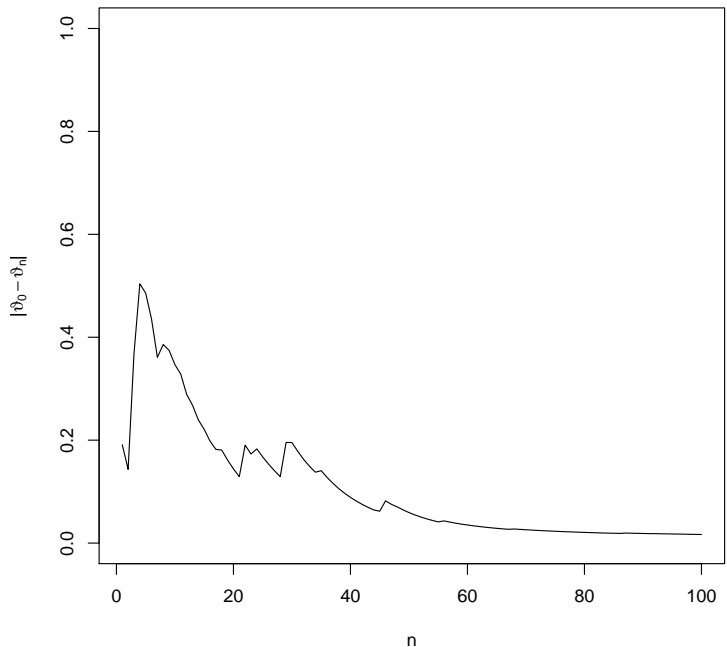}
\caption{Evolution of error  $\|\tilde\vartheta _n-\vartheta_0\|$}
\end{figure}

As can be seen the distance 
between $\vartheta_0$ and the BE after initial fluctuations quickly
decreases towards zero which illustrates the  consistency of the BE.

We also made simulations for the  MLE  $\hat\vartheta _n$
of the same parameter $\vartheta_0$. 

In what follows we present the graphs of the corresponding error obtained for the same
simulation model with $n$ running in the range $[1,100]$.


\begin{figure}[ht]
\hspace{3cm}\includegraphics[width=9cm]   {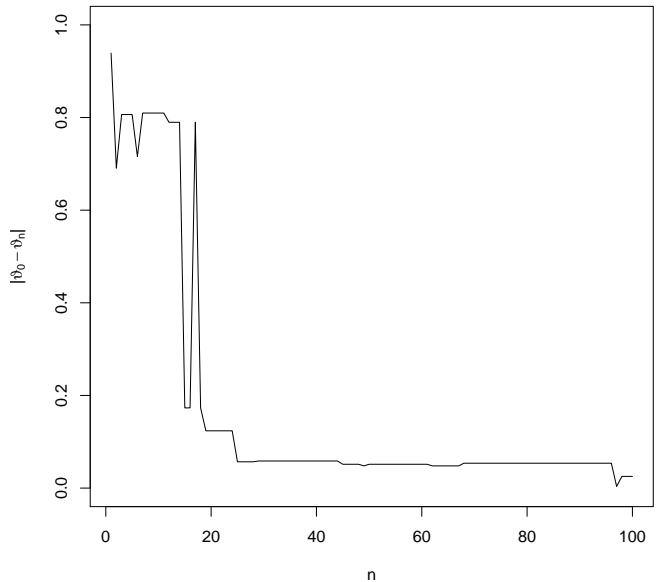}
\caption{Evolution of error  $\|\hat\vartheta _n-\vartheta_0\|$}
\end{figure}

As can be seen the
Euclidean distance between the MLE and $\vartheta_0$ quickly decreases towards
zero after initial fluctuations:

We can see that the fluctuations of the MLE at the beginning are more
important than those of the BE. 

 \section{Discussions}

Let us mention now some problems related to the study of the MLE. The main technical  difficulty to apply the
Ibragimov and Khasminskii approach in the study of the MLE in this change-point
statement is  in the checking of the tightness of the family of measures
induced by the likelihood ratio random field $Z_n\left(u\right),u\in\UU_n$ in
the space of its realizations. Recall that this is the space of surfaces with
discontinuities along some curves. 

Here we supposed that the signal and noise are of
the same magnitude $n$, where $n\rightarrow \infty $. However, in some cases the
signal can be much larger than the noise, say, 
\begin{align*}
\lambda _{j,n}\left(\vartheta ,t\right)=n\lambda \left(t-\tau
_j\left(\vartheta \right)\right)\1_{\left\{t>\tau _j\left(\vartheta
  \right)\right\}}+\lambda _0,\qquad 0\leq t\leq T,\qquad j=1,2,3.
\end{align*}
This case could  be studied as well by means of the presented method but
the limit  $Z\left(u\right),u\in\RR^2$ of the normalized     likelihood ratio
function 
$Z_n\left(u\right),u\in\UU_n$ will be different.

As mentioned in the Introduction there are several other statements related
to the problem of Poisson source localization depending on the
regularity of the signals. The cases of smooth signals and cusp-type signals
are considered in the works \cite{CK18} and \cite{CDK18} respectively. In
particular, in  \cite{CK18} the estimation of the
parameter $\vartheta _0$ by $k\geq 3$ sensors was made in two steps. First we estimate
the moments of the arrival times of the signals, say, $\bar\tau _{1,n},\ldots,\bar\tau
_{k,n}$, then given these estimators the localization $\bar\vartheta _n$ of
the source is found by solving the system of equations
\begin{align*}
\bar\tau _{1,n}^2\nu ^2=\left\|\vartheta _1-\bar\vartheta _n \right\|^2,\qquad 
\ldots,\qquad 
\bar\tau _{k,n}^2\nu ^2=\left\|\vartheta _k-\bar\vartheta _n \right\|^2.
\end{align*}
It is shown that the estimator $ \bar\vartheta _n$ is consistent and
asymptotically normal. It will be interesting to study the similar estimator
in the change-point case.

Another question concerns the robustness of the estimators (MLE and BE) with
respect to the knowledge of the model. 
Suppose that the signal $\lambda \left(t\right), t\geq 0$ is not exactly known
and we use just a constant value $\lambda _1>0$.   We can  see what are the
limits of the MLE and BE in such situations. It is known that in this case 
both estimators converge to the value $\hat\vartheta $ which minimizes the
corresponding Kulback-Leibler distance. The  one-dimensional case 
was studied in \cite{DFK01}, where it was shown that for a wide
range of values of $\lambda _1$ the BE is consistent even for the wrong
model. We could suppose that the model considered in the present work has a similar
property. Then the consistent estimation is possible in the case of
misspecification as well.

Of course a similar problem could be studied for the models of signals in
white Gaussian noise. Indeed, suppose that we have the same positions of the
source and the detectors (see Fig. 1), but the signals are Gaussian
\begin{align*}
{\rm d}X_{j,t}=S\left(t-\tau _j \left(\vartheta
\right)\right)\1_{\left\{t\geq \tau _j  \left(\vartheta
  \right)\right\}}{\rm d}t+\varepsilon {\rm d}W_{j,t}, \quad X_0=0,\quad 0\leq
t\leq T.
\end{align*}
Here $j= 1,2,3$ and   $W_{j,t}, 0\leq t\leq T, j=1,2,3$ are independent Wiener
processes. Then we can describe the properties of the MLE and BE of the
coordinates of the source in the asymptotics of small noise ($\varepsilon
\rightarrow 0$) in the cases of different regularity of the signals (see
e.g. \cite{DKKN18}). 

{\bf Acknowledgment.}  This work
was done under partial financial support of the grant of RSF 14-49-00079.

\end{document}